\documentclass[12pt]{article}
\usepackage{amsmath}
\usepackage{graphicx}
\usepackage{enumerate}
\usepackage{natbib}
\usepackage{mathrsfs}
\usepackage{amsmath}
\usepackage{amssymb,color}
\usepackage{graphicx}
\usepackage{bm}
\usepackage{epstopdf}
\usepackage{url} 
\usepackage{hyperref}
\usepackage{caption,subcaption}
\usepackage{lineno,hyperref}

\usepackage{mathrsfs}

\usepackage[nofiglist,notablist]{endfloat}

\newcommand{\bZ}{{\bm Z}}

\newcommand{\bh}{{\bm h}}

\newcommand{\ba}{{\bm a}}
\newcommand{\bb}{{\bm b}}

\newcommand{\bbeta}{\boldsymbol{\beta}}

\begin{document}

\def\spacingset#1{\renewcommand{\baselinestretch}%
{#1}\small\normalsize} \spacingset{1}
\date{}

  \title{Online survival analysis with  quantile regression}

\author{\small  Yi Deng\textsuperscript{1}~~~Shuwei Li\textsuperscript{2}~~~Liuquan Sun\textsuperscript{3}~~~ and ~~~ Baoxue Zhang\textsuperscript{1}
\\
{\small \textsuperscript{1}School of Statistics, Capital University
	of Economics and Business, Beijing, China}\\
{\small \textsuperscript{2}School of Economics and Statistics, Guangzhou University, Guangzhou, China} \\
\small \textsuperscript{3}Academy of Mathematics and Systems Science, Chinese Academy of Sciences,  \\
{\small and School of Mathematical Sciences, University of Chinese Academy of Sciences, Beijing, China}\\
}


\maketitle

\begin{abstract}
We propose an online inference method for censored quantile regression with streaming  data sets.
A key  strategy     is to approximate the martingale-based unsmooth objective function
with a quadratic loss function involving a well-justified second-order expansion.
This   enables us to derive a new online convex function  based on the current data batch and summary statistics of historical data,
thereby  achieving online updating and  occupying low storage space.
To estimate the regression parameters,
we design a novel majorize-minimize algorithm by reasonably constructing a  quadratic  surrogate objective function,  which renders a closed-form parameter update and thus reduces the computational burden notably.
Theoretically, compared to the oracle estimators  derived from analyzing
the entire raw data once, we   posit a weaker assumption on the quantile grid size and show that  the proposed online estimators  can maintain the same convergence rate and statistical efficiency.
Simulation studies and an application
demonstrate the satisfactory empirical performance and practical utilities of the proposed online method.
\end{abstract}

\noindent%
\vfill


\newpage
\spacingset{1.5} 
\section{Introduction}
Quantile regression has attracted considerable attentions  since it can offer
a comprehensive and flexible   assessment of heterogeneous covariate effects
on the response \citep{Koenker1978}.
Due to its remarkable superiorities,
numerous methods have been developed to
fit  quantile regression model  to survival data under various settings
\citep{Powell1986, Ying1995, Portnoy2003, Peng2008, Peng2009, WangWang2009,  Leng2013, LiPeng2015, Gorfine2017,  2023Wang, Chu2024}.
In particular,
under the traditional right censored data,
\cite{Portnoy2003}  proposed a redistribution-of-mass approach to estimate the censored quantile regression model.
\cite{Peng2008} developed a  martingale-based estimating equation coupled  with a $L_1$-type minimization procedure.
\cite{WangWang2009}
investigated a locally weighted estimating  equation method without requiring the global linear assumption.
\cite{Peng2009} suggested an   inverse
probability weighting estimation approach for    quantile regression  under competing risks.
\cite{LiPeng2015} proposed a copula-aided quantile modeling method to accommodate dependent censoring in semicompeting risks.


Thanks to technological
advances, large-scale streaming data are collected  in high frequency (i.e., daily, weekly or monthly)  and ubiquitous in various scientific areas,
including but not limited to financial trading \citep{Das2018}, environment assessment  \citep{2022Jiang}
and  business analysis  \citep{2024Chen}.
When large-scale data arrive in streams, conventional statistical tools
become computationally cumbersome or even infeasible due to the limited computing memory for storing the entire raw data.
Such constraints   have continuously brought out vast  online updating methods
for carrying out regression analysis   \citep{2016Schifano, 2020Xue, 2020Luo,  2021Wang, 2021Wu, 2022Jiang, 2023Sun, 2024Chen}.
For example, \cite{2016Schifano} derived an online updating method for the classical parametric  linear
regression model.
\cite{2024Yang}  proposed an online smooth backfitting method for the nonparametric generalized additive model.
\cite{2024Chen} considered a renewable quantile regression model.
For streaming survival data,  \cite{2020Xue} proposed a cumulatively updated estimating equation
 for the proportional hazards (PH) model and an online   test statistic to verify the PH assumption.
\cite{2021Wu} provided online updating methods to estimate both regression coefficients and   baseline hazard of the PH model with time-dependent covaraites.
While various analytical tools for streaming data sets have been
studied,
developing an online updating method for censored quantile regression,
a vital  model that can delineate heterogeneous covariate effects  as introduced above,
remains  a complex and  unexplored task.


This work offers an online  inference method for censored quantile regression  with streaming  data sets.
We mainly leverage  the martingale-based estimating equation proposed by \cite{Peng2008},
which involves an unsmooth objective function.
To approximate an unsmooth  function,  \cite{2022Jiang},  \cite{2023Sun} and  \cite{2023Wang}  utilized kernels that involve  selecting the optimal bandwidth. Such kernel methods  may lead to  the  substantial computational burden and estimation instability at two  tails of the kernel function.
As an alternative, we propose to approximate  the unsmooth objective function  related to the historical data  with a rigorously justified second-order expansion.
This renders a new online convex function  that only uses the current data batch and summary statistics of historical data,
thereby  achieving online updating and occupying low storage space.

There  exist  two primary challenges to construct an online  estimation method for censored quantile regression based on the work of \cite{Peng2008}.
First, since the method in \cite{Peng2008} relies on a grid-based estimation procedure,
the objective function contains the  entire raw data-based estimates.
Thus, the  entire raw data-based objective function   cannot be  directly decomposed into the summation of independent objective functions related to the separated data blocks, which is essential   to construct our online renewable estimators.
Second, since the objective function is unsmooth, it is not trivial to extract some useful summary statistics from the historical data with the traditional second-order Taylor expansion.
To tackle the first issue, by assuming that the sample size of each data block is sufficiently large, we approximate the oracle estimators derived from analyzing the entire raw data once by the  estimators related to separated  data blocks, and further verify that the entire raw data-based  objective function can be reasonably approximated by the summation of   the objective functions formed form separated  data blocks.
To solve the second issue, we  use the local quadratic approximation technique to derive a new objective function, which only involves some summary statistics including  historical data-based   weight matrix and  parameter estimators.

To estimate the regression parameters in the unsmooth objective function, we  design a majorize-minimize algorithm with a simple quadartic surrogate function, which renders a closed-form parameter update and thus reduces the computational complexity.
Compared to the offline method of \cite{Peng2008},   we    posit a weaker assumption on the quantile grid size to establish the asymptotic properties  of the proposed online estimators, including the consistency and
asymptotic normality.
In terms of the convergence rate and statistical efficiency, the proposed   estimators are shown to be  asymptotically equivalent to the oracle estimators derived from analyzing
the entire data once.
Moreover, in our computational complexity analysis, we show that  the proposed method significantly reduces both time and space complexity compared to traditional methods, leading to enhanced computational performance as also verified in the simulations.


The remainder of the article is organized as follows.
In Section 2, we
first review the classical martingale-based estimating equation for censored quantile regression and then
present  our proposed online objective function.
Section 3 provides specific algorithms used to implement the proposed online method  and discusses its computational complexity.
In Section 4, we  establish the large sample properties of the resulting estimators  and then provide a variance estimation procedure.
Simulations are conducted in Section 5 to assess  the empirical  performance of the proposed  method,
followed by a real-world data analysis in Section 6.
Section 7 provides some concluding remarks.

\section{Methodology}
\subsection{Review of martingale-based estimation}
\label{Review}

Let $T$ and $C$ be the event time of interest and right censoring time, respectively.
Define $X=\min(T, C)$ and $\Delta =I(T \le C)$, where $I(\cdot)$ is the indicator function.
Let $\bZ = (1, \widetilde{\bZ}^{\rm T})^{\rm T}$ be a $p \times 1$ vector,
where  $\widetilde{\bZ}$ is a column covariate vector.
The observed data, denoted by $\{X_{i}, \Delta_{i},\bm Z_{i}\}_{i=1}^{n}$, are $n$ independent and identically distributed replicates of   $(X, \Delta, \bm Z)$.
Herein, we make a  conditionally independent censoring assumption in the sense  that $C$ is   independent of $T$  given $\bm Z$.
Define the $\tau$th conditional quantile of $T$ given $\bm Z$ as  $Q_{T}(\tau \mid \bm Z)=\inf\{t: P( T\le t \mid \bm Z) \ge \tau\}$, where $\tau \in (0, 1)$.
A quantile regression model assumes

\begin{equation}
Q_{T}(\tau \mid \bm Z)=\exp\{\bm Z^{\rm T} \bm \beta (\tau)\},
\label{m1}
\end{equation}
where $\bbeta(\tau)$ is a $p \times 1$ vector of unknown regression coefficients, which represents the covariate effects for the $\tau$th quantile of $T$.

To estimate $\bbeta (\tau)$,   \cite{Peng2008} suggested a martingale-based estimating equation coupled with a $L_1$-type minimization procedure.
In particular,
let $N_{i}(t)=I(X_{i}\le t, \Delta_{i}=1)$ for $i=1,...,n$.
Denote by $\Lambda_{T}(t   \mid \bm Z)=-\log\{1-P( T\le t \mid \bm Z)\}$ the cumulative conditional hazard function of $T$ given $\bm Z$.
Under model \eqref{m1},    \cite{Peng2008} showed that
\begin{equation}
E\Bigg[\frac{1}{\sqrt{n}}\sum_{i=1}^{n}\bm Z_i \Big(N_i(\exp\{\bm Z_i^{\rm T}\bbeta_0(\tau)\})-\Lambda_{T}[\exp\{\bm Z_{i}^{\rm T}\bbeta_0(\tau)\} \wedge X_i  \mid \bm Z_i]\Big)\Bigg]=0,
\label{EE1}
\end{equation}
where   $\bbeta_0(\cdot)$ represents the true $\bbeta (\cdot)$ in model \eqref{m1}.  Then
\cite{Peng2008} defined an estimator of $\bbeta (\tau)$, denoted
by $\widetilde{\bbeta}_n(\tau)$, as a right-continuous piecewise-constant function,
which   jumps  on a grid  $\mathcal{M}_{K} =\{0=\tau_{0}<\tau_{1}<\cdots <\tau_{K} <1\}$.
Let   $\exp\{\bm Z^{\rm T} \widetilde{\bbeta}_n (0)\}=0$  and $H(u)=-\log(1-u)$.
\cite{Peng2008}
proposed to  obtain $\widetilde{\bbeta}_n (\tau_k)  \,   (k = 1, \ldots , K)$ by sequentially   solving the following monotone estimating equation
$$
\frac{1}{n}\sum_{i=1}^{n}\bm Z_{i}\Bigg(N_{i}(\exp\{\bm Z_{i}^{\rm T} \bm \beta(\tau_k)\})-\sum_{r=0}^{k-1} I\left[X_{i}\ge \exp\{\bm Z_{i}^{\rm T}\widetilde{\bbeta}_n(\tau_{r})\}\right] \{H(\tau_{r+1}) - H(\tau_{r})\} \Bigg)=0.
$$
Equivalently,  for each  $k = 1, \ldots , K$, $\widetilde{\bbeta}_n (\tau_k)$ can be obtained by minimizing the following convex objective function
\begin{align}
l_{kn}(\bbeta(\tau_k)) =&\frac{1}{n}\sum_{i=1}^{n}\Bigg[\Delta_{i}|\log X_{i}-\bbeta(\tau_k)^{\rm T}\bm Z_{i}|+\bbeta(\tau_k)^{\rm T}\Delta_{i} \bm  Z_{i} \Bigg.  \notag \\ &
\Bigg.-2\bbeta(\tau_k)^{\rm T}\bm  Z_{i}\Bigg(\sum_{r=0}^{k-1} I\left[ X_{i}\ge \exp\{\bm  Z_{i}^{\rm T}\widetilde{\bbeta}_n(\tau_r)\}\right] \{H(\tau_{r+1}) - H(\tau_{r})\}\Bigg)\Bigg].\label{loss1}
\end{align}
Such a minimization can be readily accomplished with the Barrodale-Roberts algorithm \citep{1974Barrodale}, which is available in standard statistical software,  such as the $rq(\cdot)$ function in R package $quantreg$.
By assuming that the quantile grid size $m=\max \{\tau_k-\tau_{k-1}; k=1, \ldots ,K\}=o(n^{-1/2})$,
\cite{Peng2008} showed that
the  estimator $\widetilde{\bbeta}_{n}(\tau)$ with $\tau \in (0,  \tau_K]$ is
$\sqrt{n}$-consistent.

\subsection{Proposed online estimation}
\label{ProMethod}

In large-scale streaming data setting,
data often arrive in chunks  at  high frequency (i.e., daily, weekly or monthly).
The estimation procedure proposed by \cite{Peng2008}   requires
accessing all the raw data directly and   high storage space.
In order to circumvent  this obstacle,
we herein propose an online updating method that is based on the current data batch and summary statistics of historical data,
thereby    warranting low memory and computational complexity.
For $b=1,2, \ldots, B$,
let $D_b$ denote all the collected data in the $b$th  batch,  $n_b$ be the sample size of $D_b$, and $N_b=\sum_{j=1}^{b} n_j$ be the cumulative sample size until the accumulation point $b$.
For each  $ k = 1, \ldots , K$,
 $\widetilde{\bbeta}_{n}(\tau_k)$  denotes the oracle estimator based on all the data with the sample size $n = N_B$, which is obtained by minimizing \eqref{loss1}.

To derive an online  estimator for   model \eqref{m1}, we begin
with considering a simple scenario that involves only two data batches   $D_1$ and $D_2$, where $D_2$ arrives after $D_1$, that is, $B=2$.
For each  $ k = 1, \ldots , K$, let $\widehat{\bbeta}_{N_b} (\tau_k)$  be the online renewable estimator at cumulative point $b=1,2, \ldots, B$. $\widehat{\bbeta}_{N_b} (0)$ satisfies $\exp\{\bm Z^{\rm T} \widehat{\bbeta}_{N_b} (0)\}=0$.
Based on the first data batch, we can obtain the online estimator   $\widehat{\bbeta}_{N_1}(\tau_k)=\widetilde{\bbeta}_{n_1}(\tau_k)={\rm argmin}_{\bbeta(\tau_k) \in \mathcal{R}^p} l_{k n_1}(\bbeta(\tau_k))$ for each $k=1, \ldots, K$, where $l_{k n_1}(\bbeta(\tau_k))$ is defined in \eqref{loss1}.
When $D_2$ arrives,
the oracle estimator based on the entire raw data $D_1  \cup D_2$ is
$\widetilde{\bbeta}_{N_2}(\tau_k)
=  {\rm argmin}_{\bbeta(\tau_k) \in \mathcal{R}^p} l_{k N_2}(\bbeta(\tau_k))$ for each $k=1, \ldots, K$.
Let $\widetilde{\bbeta}_{kN_B}=\{\widetilde{\bbeta}_{N_B}(\tau_r), 0< r \le k-1\}$ if no confusion will be caused.
Notably,  $l_{k N_2}(\bbeta(\tau_k))$      can be naturally decomposed  as
$$
\frac{1}{N_2}\big\{n_1 \widetilde{l}_{kn_1}(\bbeta(\tau_k), \widetilde{\bbeta}_{kN_2})+n_2 \widetilde{l}_{kn_2}(\bbeta(\tau_k), \widetilde{\bbeta}_{kN_2})\big\},
$$
where
$$
\begin{aligned}
\widetilde{l}_{kn_b}(\bbeta(\tau_k), \widetilde{\bbeta}_{kN_2}) =&\frac{1}{n_b}\sum_{i=1}^{n_b}\Bigg[\Delta_{i}|\log X_{i}-\bbeta(\tau_k)^{\rm T}\bm Z_{i}|+\bbeta(\tau_k)^{\rm T}\Delta_{i} \bm  Z_{i}\Bigg. \\ &
\Bigg.-2\bbeta(\tau_k)^{\rm T}\bm  Z_{i}\Bigg(\sum_{r=0}^{k-1}I\Big[X_{i}\ge \exp\{\bm  Z_{i}^{\rm T}\widetilde{\bbeta}_{N_2}(\tau_r)\}\Big]\{H(\tau_{r+1}) - H(\tau_{r})\}\Bigg)\Bigg],
\end{aligned}
$$
for $b=1$ and 2.
However, since $l_{k N_2}(\bbeta(\tau_k))$
involves $\widetilde{\bbeta}_{N_2}(\tau_r)$ $(0< r \le k-1)$, the estimator of $\bbeta(\tau_r)$
 based on the entire raw data  $D_1 \cup D_2$, $l_{kN_2}(\bbeta(\tau_k))$  cannot be directly written as  $\{n_1 l_{kn_1}(\bbeta(\tau_k))+n_2 l_{kn_2}(\bbeta(\tau_k))\}/N_2$, which is essential to the development of our online estimation method.
To address this issue, we show in  the proposition below  that
$\widetilde{l}_{kn_1}(\bbeta(\tau_k), \widetilde{\bbeta}_{kN_2})$ and $\widetilde{l}_{kn_2}(\bbeta(\tau_k), \widetilde{\bbeta}_{kN_2})$
 can  be  approximated
  by   $l_{kn_1}(\bbeta(\tau_k))$ and  $l_{kn_2}(\bbeta(\tau_k))$, respectively.
 Let $\| \cdot \|$ denote the Euclidean  norm.\\

 \noindent{\bf Proposition 1}. If $n_b \rightarrow \infty$ for $b=1, \ldots, B$,  and there exists a constant $c>0$ such that $\|\bbeta(\tau_k)\| \le c(N_B^{-1}\sum_{b=1}^{B} \sqrt{n_{b}})^{-1}$ for each $k=1, \ldots, K$, then
$l_{k N_B}(\bbeta(\tau_k))$   can be approximated by $N_B^{-1} \sum_{b=1}^{B} n_b$ $l_{k n_b}(\bbeta(\tau_k))$
with  the quantile grid size $m=\max \{\tau_k-\tau_{k-1}; k=1, \ldots ,K\}   =o(\tilde{n}^{-1/2})$, where  $\tilde{n}=\max \{n_b; b=1, \ldots, B\}$. \\

The proof of Proposition 1  is sketched in Section A of the Supplementary Material.
This proposition states that $l_{k N_2}(\bbeta(\tau_k))$ can be decomposed as
a weighted average of $l_{k n_1}(\bbeta(\tau_k))$   and $l_{k n_2}(\bbeta(\tau_k))$ as follows
\begin{equation}
\frac{1}{N_2}\{n_1 l_{kn_1}(\bbeta(\tau_k)) +n_2 l_{kn_2}(\bbeta(\tau_k))\}.
\label{loss2}
\end{equation}

 Moreover, in  the proposed  online framework,
we only posit that the quantile grid size $m = o(\tilde{n}^{-1/2})$,    whereas
\cite{Peng2008}'s oracle method based on the entire raw data
requires a more stringent assumption  that  $\max \{\tau_k-\tau_{k-1}; k=1, \ldots ,K\}=o(N_B^{-1/2})$  since $\tilde{n}  <<  N_B$ as $B \rightarrow \infty$ in the online   setting,
where the symbol $<<$ denotes ``much less than''.

To update each $\bbeta (\tau_k)$ in an online manner, when $D_2$ arrives,
we hope to utilize $D_2$ and some historic summary statistics based on $D_1$ instead of the raw data.
Note that the minimization of      $\frac{1}{N_2}\{n_1 l_{kn_1}(\bbeta(\tau_k))+n_2 l_{kn_2}(\bbeta(\tau_k))\}$ with respect to $\bbeta(\tau_k)$ is equivalent to minimizing
\begin{equation*}
\begin{aligned}
\frac{n_1}{N_2}\{l_{kn_1}(\bbeta(\tau_k))-l_{kn_1}(\widehat{\bm \beta}_{N_1}(\tau_k))\}+\frac{n_2}{N_2}l_{kn_2}(\bbeta(\tau_k)).
\end{aligned}
\end{equation*}
In Proposition 2 below,
we show that $l_{kn_1}(\bbeta(\tau_k))-l_{kn_1}(\widehat{\bm \beta}_{N_1}(\tau_k))$
can be approximated by a quadratic loss function,
which only involves summary statistics based on $D_1$ instead of the raw data.
Let $\widehat{\bbeta}_{N_b} (\tau)$ be a right-continuous piecewise-constant function that jumps only on a grid  $\mathcal{M}_{K}$ at cumulative point $b=1,2, \ldots, B$.
\\

\noindent{\bf Proposition 2.} If $n_1 \rightarrow \infty$ and  $\|\bbeta(\tau)-\bbeta_{0}(\tau)\|$ is bounded for $\tau \in (0, \tau_K]$,
then we have
$$
l_{kn_1}(\bbeta(\tau))-l_{kn_1}(\widehat{\bbeta}_{N_1}(\tau))
=\{\bbeta(\tau)-\widehat{\bbeta}_{N_1}(\tau)\}^{\rm T} \Gamma_{n_1}(\tau)\{\bbeta(\tau)-\widehat{\bbeta}_{N_1}(\tau)\}+o_p(1),
$$
where $\Gamma_{n_b}(\tau)=n_b^{-1}\sum_{i=1}^{n_b} E[\bm Z_{i}^{(b)} \bm Z_{i}^{(b)\rm T} f(\exp\{\bm Z^{(b)\rm T}\bm \beta_{0}(\tau) \}|\bm Z_{i}^{(b)} )\exp\{\bm Z^{(b)\rm T}\bm \beta_{0}(\tau) \}]$,  $\bm Z_{i}^{(b)}$ is the covariate vector of the $i$th subject in the $b$th data batch,
$f(x \mid \bZ)={\rm d}F(x \mid \bZ)/{\rm d}x$ and $F(t \mid \bZ)=P(X \le t, \Delta=1 \mid \bZ)$.
\\

The proof of Proposition 2 is given in Section A of Supplementary Material.
Let $\widehat{\Gamma}_{n_b}(\tau)$ be the consistent estimate for $\Gamma_{n_b}(\tau)$. $\widehat{\Gamma}_{n_b}(\tau)$ is also a right-continuous piecewise-constant function that jumps only on a grid  $\mathcal{M}_{K}$.
For each $k=1, \ldots, K$, according to   Proposition 2 and replacing  the weight matrix $\Gamma_{n_1}(\tau)$ with $\widehat{\Gamma}_{n_1}(\tau_k)$ for $\tau \in [\tau_k, \tau_{k+1})$,  $l_{k N_2}(\bbeta(\tau_k))$ 
can be reasonably approximated by
$$
\frac{n_1}{N_2}\{\bbeta(\tau_k) - \widehat{\bm \beta}_{N_1}(\tau_k)\}^{\rm T}\widehat{\Gamma}_{n_1}(\tau_k)\{\bbeta(\tau_k) - \widehat{\bm \beta}_{N_1}(\tau_k)\}+\frac{n_2}{N_2}l_{kn_2}(\bbeta(\tau_k))
$$
with quantile grid size $m = o(\tilde{n}^{-1/2})$.
Therefore, we can update   $\bbeta (\tau_k)$ in an online manner by leveraging the current data   $D_2$ and   summary statistics $\{n_1, \widehat{\bbeta}_{N_1}(\tau_k), \widehat{\Gamma}_{n_1}(\tau_k)\}$ related to the data  $D_1$ for each $k=1, \ldots, K$.

Similarly, let $\widehat{\Gamma}_{N_B}(\tau)$, a right-continuous piecewise-constant function that jumps only on a grid  $\mathcal{M}_{K}$, be the consistent estimator for $\Gamma_{N_B}(\tau)$ with
$$
\Gamma_{N_B}(\tau)=\frac{1}{N_B}\sum_{b=1}^{B}\sum_{i=1}^{n_b} E[\bm Z_{i}^{(b)} \bm Z_{i}^{(b)\rm T}  f(\exp\{\bm Z^{(b)\rm T}\bm \beta_{0}(\tau) \} \mid \bm Z_{i}^{(b)} )\exp\{\bm Z^{(b)\rm T}\bm \beta_{0}(\tau) \}].
$$
Considering $B$ data batches in total,
the proposed online estimator of  $\bbeta (\tau_k)$, denoted by
$\widehat{\bbeta}_{N_B}(\tau_k)$, is the minimizer of the following loss function
\begin{align}
    G_{kB}(\bbeta(\tau_k))=& \frac{N_{B-1}}{N_B}\{\bbeta(\tau_k)-\widehat{\bm\beta}_{N_{B-1}}(\tau_k)\}^{\rm T}\widehat{\Gamma}_{N_{B-1}}(\tau_k)\{\bbeta(\tau_k)-\widehat{\bm\beta}_{N_{B-1}}(\tau_k)\} \notag \\&+\frac{n_B}{N_B}l_{k n_B}(\bbeta(\tau_k)) \label{lossG}
\end{align}
for each $k=1, \ldots, K$.
Thus, the proposed online estimation procedure depends only on current data  $D_{B}$ and summary statistics $\{N_{B-1}, \widehat{\bbeta}_{N_{B-1}}(\tau_k), \widehat{\Gamma}_{N_{B-1}}(\tau_k)\}$ related to  the data  at cumulative point $B-1$ for $k=1, \ldots, K$.
Differentiating the loss function \eqref{lossG} with respect to $\bbeta(\tau_k)$ for each $k=1, \ldots, K$,
the online estimator $\widehat{\bm \beta}_{N_B}(\tau_k)$  is the solution of the following estimating equation
$$
H_{kB}(\bbeta(\tau_k))=\frac{N_{B-1}}{N_B}\widehat{\Gamma}_{N_{B-1}}(\tau_k)\{\bbeta(\tau_k)-\widehat{\bbeta}_{N_{B-1}}(\tau_k)\}+\frac{n_B}{N_B}S_{k n_B}(\bbeta(\tau_k))=0,
$$
where
$$
S_{k n}(\bbeta(\tau_k))=\frac{1}{n}\sum_{i=1}^{n}\bm Z_{i}\Bigg(N_{i}(\exp\{\bm Z_{i}^{\rm T} \bbeta(\tau_k)\})-\sum_{r=0}^{k-1} I\left[X_{i}\ge \exp\{\bm Z_{i}^{\rm T}\widetilde{\bbeta}_n(\tau_{r})\}\right] \{H(\tau_{r+1}) - H(\tau_{r})\} \Bigg).
$$

\section{Computational Methods and  Complexity Analysis}
\subsection{Calculation of   $\widehat{\Gamma}_{N_B}(\tau)$}
\label{ReLS}

Note that $\Gamma_{N_B}(\tau)$ contains the unknown function $f(t \mid \bm Z)$,   the conditional density function of $T$  given $\bm Z$.
Thus, to obtain  $\widehat{\Gamma}_{N_B}(\tau_k)$ $(k=1, \ldots, K)$,
we first need to estimate $f(t \mid \bm Z)$.
Traditional nonparametric kernel methods used in \citep{WangWang2009} and others are cumbersome and unstable since it  depends on the intensive
bandwidth selection and suffers from the curse of dimensionality.
Herein,  motivated by the work of \citep{2008Zeng}, we propose to directly calculate $\widehat{\Gamma}_{N_B}(\tau_k)$
with   an efficient and
stable resampling least squares (ReLS) strategy, which essentially bypasses
estimating $f(t \mid \bm Z)$.

In particular,
according to Lemma S1 in Section C of the Supplementary Material, when $B=1$, we have
$$
\sqrt{N_1}S_{kN_1}(\widehat{\bbeta}_{N_1}(\tau))=\sqrt{N_1}S_{kN_1}(\bbeta_{0}(\tau))+\Gamma_{n_1}(\tau)\sqrt{N_1}\{\widehat{\bbeta}_{N_1}(\tau)-\bbeta_0(\tau)\}+o_p(1),
$$
and for $B>1$,
$$
\sqrt{N_B}H_{kB}(\widehat{\bbeta}_{N_B}(\tau))=\sqrt{N_B}H_{kB}(\bbeta_{0}(\tau))+\Gamma_{N_B}(\tau)\sqrt{N_B}\{\widehat{\bbeta}_{N_B}(\tau)-\bbeta_0(\tau)\}+o_p(1),
$$
where $\tau \in (0, \tau_K]$.
Thus, for each $k=1, \ldots, K$,  if $B=1$,
$$
\sqrt{N_1}S_{kN_1}(\widehat{\bm \beta}_{N_1}(\tau_k)+\bm\xi_k/\sqrt{N_{1}})-\sqrt{N_1}S_{kN_1}(\widehat{\bm \beta}_{N_1}(\tau_k))=\Gamma_{n_1}(\tau_k)\bm\xi_k+o_p(1),
$$
and if $B>1$,
$$
\sqrt{N_B}H_{kB}(\widehat{\bm \beta}_{N_B}(\tau_k)+\bm\xi_k/\sqrt{N_{B}})-\sqrt{N_b}H_{kB}(\widehat{\bm \beta}_{N_B}(\tau_k))=\Gamma_{N_B}(\tau_k)\bm\xi_k+o_p(1),
$$
where $\bm \xi_k$ is a $p$-dimensional zero-mean random vector independent of the steaming data.
Since $S_{kN_1}(\widehat{\bbeta}_{N_1}(\tau_k))=0$ and $H_{kB}(\widehat{\bbeta}_{N_B}(\tau_k))=0$, we have $$\sqrt{N_1}S_{kN_1}(\widehat{\bm \beta}_{N_1}(\tau_k)+\bm\xi_k/\sqrt{N_{1}})=\Gamma_{n_1}(\tau_k)\bm\xi_k+o_p(1),$$ and $$\sqrt{N_B}H_{kB}(\widehat{\bm \beta}_{N_B}(\tau_k)+\bm\xi_k/\sqrt{N_{B}})=\Gamma_{N_B}(\tau_k)\bm\xi_k+o_p(1).$$
Regarding each row of the matrix $\Gamma_{N_B}(\tau_k)$ as the unknown parameter vector of a linear function and $\sqrt{N_1}S_{kN_1}(\widehat{\bm \beta}_{N_1}(\tau_k)+\bm\xi_k/\sqrt{N_{1}})$ or $\sqrt{N_B}H_{kB}(\widehat{\bm \beta}_{N_B}(\tau_k)+\bm\xi_k/\sqrt{N_{B}})$ as the response, we can obtain the estimate of $\Gamma_{N_B}(\tau_k)$ as follows
$$
\widehat{\Gamma}_{N_B}(\tau_k)=(\bm \xi_k^{\rm T}\bm \xi_k)^{-1}\bm \xi_k^{\rm T}W_B(\tau_k),
$$
where $W_1(\tau_k)=\sqrt{N_1}S_{kN_1}(\widehat{\bbeta}_{N_1}(\tau_k)+\bm\xi_k/\sqrt{N_{1}})$ if $B=1$ and $W_B(\tau_k)=\sqrt{N_B}H_{kB}(\widehat{\bm \beta}_{N_B}(\tau_k)+\bm\xi_k/\sqrt{N_{B}})$ if $B>1$.
Taking for example the setting of $B>1$, we can obtain $\widehat{\Gamma}_{N_B}(\tau_k)$ by the following steps.
First, we generate $s$ independent realizations $\{\bm\xi_k^{(1)}, \ldots,\bm\xi_k^{(s)}\}$ of
the random vector $\bm\xi_k$, where $s$ is a large positive integer.
Next, obtain $\sqrt{N_B}H_{kB}(\widehat{\bm \beta}_{N_B}(\tau_k)+\bm\xi_k^{(\phi)}/\sqrt{N_{B}})$  for $\phi=1, \ldots, s$.
Finally, we calculate the least squares estimate
of $[\sqrt{N_B}H_{kB}( \widehat{\bm \beta}_{N_B}(\tau_k)
+\bm\xi_k^{(\phi)}/\sqrt{N_{B}})]_l$,
where $[\ba]_l$ denotes the $l$th element of a $p$-dimensional column vector
$\ba$    for $l = 1, \ldots, p$.
Let the $l$th row of $\widehat{\Gamma}_{N_B}(\tau_k)$,
denoted by  $[\widehat{\Gamma}_{N_B}(\tau_k)]_l$, be the $l$th least squares
estimate derived above.
Then we have $[\widehat{\Gamma}_{N_B}(\tau_k)]_l=(\bm \xi_k^{\rm T}\bm \xi_k)^{-1}\bm \xi_k^{\rm T}[W_B(\tau_k)]_l$,
where  $\bm\xi_k=[\bm\xi_k^{(1)}, \ldots ,\bm\xi_k^{(s)}]^{\rm T}$ and $[W_B(\tau_k)]_{l}=([\sqrt{N_B}H_{kB}(\widehat{\bm \beta}_{N_B}(\tau_k)+\bm\xi_k^{(1)}/\sqrt{N_{B}})]_l, \ldots, [\sqrt{N_B}H_{kB}(\widehat{\bm \beta}_{N_B}(\tau_k)+\bm\xi_k^{(s)}/\sqrt{N_{B}})]_l)^{\rm T}$.
The  consistency of $\widehat{\Gamma}_{N_B}(\tau)$ with $\tau \in (0, \tau_K]$
was established in Lemma S2 in Section C of the Supplementary Material.

\subsection{Calculation of  the proposed online estimator $\hat{\bbeta}_{N_B}(\tau)$}
Based on the arguments in Section
\ref{ProMethod}, we can obtain $\widehat{\bm \beta}_{N_B}(\tau_k)$ $(k=1, \ldots, K)$ by minimizing $G_{kB}(\bbeta(\tau_k))$ in \eqref{lossG} sequentially, that is
$$
\begin{aligned}
\widehat{\bbeta}_{N_B}(\tau_k)={\rm argmin}_{\bbeta(\tau_k) \in \mathcal{R}^{p}}\Bigg\{&\frac{N_{B-1}}{N_B}\{\bbeta(\tau_k)-\widehat{\bbeta}_{N_{B-1}}(\tau_k)\}^{\rm T}\widehat{\Gamma}_{N_{B-1}}(\tau_k)\{\bbeta(\tau_k)-\widehat{\bbeta}_{N_{B-1}}(\tau_k)\}
\\
&+\frac{n_B}{N_B}l_{k n_B}(\bbeta(\tau_k))\Bigg\}.
\end{aligned}
$$

 Since $G_{kB}(\bbeta(\tau_k))$  is discontinuous and difficult to minimize,
we herein propose  a  majorize-minimize (MM) algorithm to obtain each $\widehat{\bbeta}_{N_B}(\tau_k)$.
Specifically, for each $k=1, \ldots, K$, to construct a proper majorizing function for $G_{kB}(\bbeta(\tau_k))$ around $\bbeta^{(0)} (\tau_k)$, where $\bbeta^{(0)} (\tau_k)$ is the initial value of $\bbeta  (\tau_k)$, we define an isotropic
quadratic function for some $\omega>0$,
$$
\begin{aligned}
F_{kB}(\bbeta(\tau_k); \omega, \bbeta^{(0)} (\tau_k))= &G_{kB}(\bbeta^{(0)} (\tau_k))+\langle 2H_{kB}(\bbeta^{(0)} (\tau_k)), \bbeta(\tau_k)- \bbeta^{(0)} (\tau_k)\rangle \\&+\frac{\omega}{2}\|\bbeta(\tau_k)-\bbeta^{(0)} (\tau_k)\|^{2},
\end{aligned}
$$
 where $\langle \ba, \bb\rangle$ denotes the inner product of  vectors $\ba$ and   $\bb$. It is easy to see that $$F_{kB}(\bbeta^{(0)} (\tau_k); \omega, \bbeta^{(0)} (\tau_k))=G_{kB}(\bbeta^{(0)} (\tau_k)).$$
Theoretically,
the quadratic coefficient $\omega>0$ should be sufficiently large such that the local majorization property $F_{kB}(\bbeta^{(1)} (\tau_k); \omega, \bbeta^{(0)} (\tau_k))\ge G_{kB}(\bbeta^{(1) }(\tau_k))$ holds, where $\bbeta^{(t)} (\tau_k)$ denotes the update of $\bbeta (\tau_k)$ at the $t$th iteration with $t \geq 1$.
In implementing the proposed algorithm,   with a relatively small value $\omega=\omega_0$ (e.g., 0.5), we  iteratively multiply $\omega_0$  by a factor of $\lambda>1$ and compute
$\bbeta^{(1)} (\tau_k)=\bbeta^{(0)} (\tau_k)-2H_{kB}(\bbeta^{(0)} (\tau_k))/w_{\alpha}$ with $\omega_\alpha =\lambda^{\alpha}\omega_{0}$ until the local majorization property holds,
where $\alpha$ is nonnegative integer sequence started from 0 and increased by 1.
At the $t$th iteration,
we can derive the majorizing function for $G_{kB}(\bbeta(\tau_k))$ around $\bbeta^{(t)}(\tau_k)$  as
$$
\begin{aligned}
F_{kB}(\bbeta(\tau_k); \omega, \bbeta^{(t)}(\tau_k))=& G_{kB}(\bbeta^{(t)}(\tau_k))+\langle 2H_{kB}(\bbeta^{(t)}(\tau_k)), \bbeta(\tau_k)- \bbeta^{(t)}(\tau_k)\rangle \\&+\frac{\omega}{2}\|\bbeta(\tau_k)-\bbeta^{(t)}(\tau_k)\|^{2},
\end{aligned}
$$
and obtain $\bbeta^{(t+1)}(\tau_k)$ as above.
Repeating this procedure yields a sequence of
iterates $\{\bbeta^{(t)}(\tau_k)\},t=0,1,\ldots$ until the stopping criterion is met,
 say $\|\bbeta^{(t+1)}(\tau_k)-\bbeta^{(t)}(\tau_k)\|\le 10^{-4}$.


\subsection{Computational complexity analysis}

Computational complexity is a crucial metric for evaluating  an algorithm's
capacity in an online framework, primarily consisting of the required storage space and computation  time complexity.
We first quantify the required storage spaces for our proposed online method and the oracle method based on the   entire raw data.
In particular, the oracle method needs to store $O(N_B)$ data
while our proposed method only occupies $O(n_B+2K p^{2}+ K p))$ volume.
This is because   our proposed method only stores the last data batch and   summary statistics $\{\widehat{\bbeta}_{N_{B-1}}(\tau_k), \widehat{\Gamma}_{N_{B-1}}(\tau_k),\sum_{b=1}^{B-1}\widetilde{\Gamma}_{n_b}(\tau_k)V_{n_b}(\tau_k)\widetilde{\Gamma}_{n_b}(\tau_k)\}$ for $k=1, \ldots, K$.
More specifically, consider the setting of  $m =o(\tilde{n}^{-1/2})$
and finite $p$,  the space complexity of our online method is $O(n_B+K)$ with $K=\tau_K / m$. The oracle method will eventually become unworkable  since $M << N_B$ as $B \rightarrow \infty$, where $M$  is the computer memory.  This states that the growing data blocks demand memory that is far beyond the computer's capacity.
At the cumulative point $B$, the available memory space and the empirical relative space loaded for running the algorithm of the proposed method are $(M-n_B)$ and $(2Kp^2+Kp)/(M-n_B)$ respectively,  where $n_B$  and $K$ are negligible compared to $M$ in practice.
Note that it is infeasible to keep the oracle quantile grid size  assumption $o(N_B^{-1/2})$, which implies that we need storing and estimating at least $O(N_B^{1/2})$ estimators, leading to inadequate computing resources.
For the proposed method with a weaker assumption $m=o(\tilde{n}^{-1/2})$,  we
release  a large amount of memory and make our estimation procedure workable since $K<<N_B$ and $\tilde{n}<<N_B$ as $B \rightarrow \infty$.

Computational time complexity depends on the time cost at each iteration
of the used optimization algorithm  and the iteration number  to achieve convergence.
Since our online method only involves inverting low dimensional   matrices, the time cost is mainly dominated by the sample size, the calculation of the summary statistics and the quantile grid size.
We quantify the computational time cost at the accumulation point $B$.
The time complexity of the oracle model is about $O(\widetilde{K} (L+1)\{N_B p^2 + p^3\})$ since the Barrodale-Roberts algorithm employed in the oracle method involves $O(N_B p^2+p^3)$ calculations  and the estimation of covariance matrix entails the resampling procedure with $L$ repetitions at a specific quantile level, where $\widetilde{K}$ is the number of the quantile points in the quantile grid set of oracle method.
By the assumption that the quantile grid size is $o(N_B^{-1/2})$  of \citep{Peng2008},  one has $\widetilde{K}=\tau_K/\widetilde{m}$ with $\widetilde{m}=o(N_B^{-1/2})$.
The time consuming of our online method mainly consists of the time cost of the proposed MM algorithm, the calculation of weight matrix $\widehat{\Gamma}_{N_{B-1}}(\tau_k)$ $(k=1, \ldots, K)$ and the covariance matrix estimation with $L$ times' resampling.
For the proposed MM algorithm, the time complexity is $O(R K n_B p)$, where $R$ denotes the required iteration number to achieve convergence.
The time complexity of the
least squares resampling used in Section  \ref{ReLS} at a given quantile level is $O(s^{2} p^2+s n_B p)$, in which the least squares method costs $O(s^2 p^2)$ and the calculation of $H_{kb}(\cdot)$ costs $O(sn_B p)$ for $s$ replicates. The bootstrap procedure involves $O(KL\{n_B p^2 + p^3\})$.
Thus, the overall time complexity of our online method is $O(R K n_B p +K \{s^{2} p^2+s n_B p\}+KL\{n_B p^2 + p^3\})$.
More specifically, in the setting of finite $p$, $R$, $s$ and $L$,
our online method reduces the time complexity to $O(K n_B)$
in contrast to the oracle method.
Since $m=o(\tilde{n}^{-1/2})$  and $\tilde{n} << N_B$ as $B \rightarrow \infty$,  we have $K << \widetilde{K}$ and thus the proposed method significantly reduces time computational complexity.

\section{Asymptotic Properties and Variance Estimation}
\subsection{Asymptotic properties}
Let $\widetilde{F}(t  \mid \bm Z)=P(X > t \mid \bm Z)$ and $\widetilde{f}(t \mid \bm Z)={\rm d}\widetilde{F}(t \mid \bm Z)/{\rm d}t$.
Define
$\mathcal{B}(\bh)=\{\bh {\in} \mathcal{R}^{p}: \inf_{\tau \in (0, \tau_K]} \|\bm \upsilon(\bh)-\bm \upsilon(\bbeta_{0}(\tau))\|\le d\}$ for $d > 0$ with $\bm \upsilon (\bh)=E[\bZ \exp\{\bZ^{\rm T} \bh\}]$, and  $\mathcal{F}=\{g: [0, \tau_K] \rightarrow \mathcal{R}^{p}, g\, \text{is} \, \text{left continuous} \, \text{with} \, \text{right} \, \text{limit, and}\, g(0)=0 \}$.
Let $A(\bh)=E[\bm Z \bm Z^{\rm T}$ \\
$ f(\exp\{\bm Z^{\rm T}\bh\} \mid \bm Z)\exp\{\bm Z^{\rm T}\bh\}]$ and $J(\bh)=E[\bm Z \bm Z^{\rm T} \widetilde{f}(\exp\{\bm Z^{\rm T}\bh\} \mid \bm Z )\exp\{\bm Z^{\rm T}\bh\}]$.
Let
$\Phi$ be a  linear operator that maps from $\mathcal{F}$ to $\mathcal{F}$ such that  $\Phi(g)(\tau)=\int_{0}^{\tau}\mathcal{I}(s, \tau) {\rm d}g(s)$ for $g \in \mathcal{F}$,
where $\mathcal{I}(s, \tau)=\prod_{u \in (s, \tau]}[I_p+J\{\bbeta_0(u)\}A(\bbeta_0(u))^{-1}{\rm d}H(u)]$ and $I_p$ is the $p \times p$ identity matrix.
To establish the asymptotic properties of the proposed online estimator $\widehat{\bbeta}_{N_B} (\tau)$, we use the following standard regularity conditions.
\begin{itemize}
\item[(C1)] The condition density function $f(t \mid \bm Z)$ and $\widetilde{f}(t \mid \bm Z)$ are bounded above 0 for any $t >0$ and $\bm Z$.
\item[(C2)] The covariate vector $\bm Z$ belongs to a compact set in $\mathcal{R}^p$,  $E(\bm Z \bm Z^{\rm T})>0$  and  $\|\bm Z\|$ is finite.

\item[(C3)] The density $f(\exp\{\bm Z^{\rm T} \bh\} \mid \bm Z) >0$ and $A(\bh)^{-1}J(\bh)$
    is uniformly bounded for any $\bh \in \mathcal{B}(\bh)$.


\item[(C4)] $\inf_{\tau \in [v, \tau_K]} eigmin A(\bbeta_0(\tau))>0$ for any $v \in (0, \tau_K)$, where $eigmin(\cdot)$ denotes the minimum eigenvalue of a matrix, and each component of $A(\bbeta_0(\tau))$ is a Lipschitz function of $\tau$.

\item[(C5)]    For each $b=1, \ldots, B$,      $n_b \rightarrow \infty$,  and  $m=\max \{\tau_k-\tau_{k-1}; k=1, \ldots ,K\}   =o(\tilde{n}^{-1/2})$.\\

\end{itemize}

Conditions (C1)--(C3) are  commonly used  in the literature of  censored quantile regression \citep{Peng2008, 2018Zheng,2022He},
which are required  to derive the consistency and asymptotic normality of an estimator based on each data batch.
These conditions  are also used   to establish the consistency of the proposed online  estimator.
The first part of condition (C4) guarantees the invertibility of $\Gamma_{N_b}(\tau)$ for $b= 1, \ldots, B$ and the identifiability of $\{\bbeta_0 (\tau), \tau \in (0, \tau_K]\}$  \citep{Peng2008}, whereas the second part is essential to establish the consistency of $\widehat{\Gamma}_{N_b}(\tau)$.
The first part of condition (C5) requires that
the sample size of each data batch  is sufficiently large,
which is useful to derive a consistent estimator related to each data batch
and the  subsequent consistent online estimator.
The second part of condition (C5) assumes that   quantile grid size $m$ is sufficiently small, which is used to show that the proposed online estimator is $\sqrt{N_B}$-consistent.
The following two theorems  describe
the asymptotic behaviors of $\widehat{\bbeta}_{N_B}(\tau)$.
\\


\noindent{\bf Theorem 1.} Under the regularity conditions (C1)--(C5),
$\sup_{\tau \in (0, \tau_K]}\|\widehat{\bbeta}_{N_B}(\tau)-\bm \beta_0(\tau)\|  \rightarrow 0$ in probability.
\\

\noindent{\bf Theorem 2.} Under the regularity conditions (C1)--(C5), $\sqrt{N_B}\{\widehat{\bbeta}_{N_B}(\tau)-\bm \beta_0(\tau)\}$ converges weakly to a mean zero Gaussian process $A(\bbeta_0(\tau))^{-1}\Phi(\Omega(\tau))$ for $\tau \in (0, \tau_K]$, where $\Omega(\tau)$ is a Gaussian process with mean $0$ and covariance $\Sigma(s, t)$ for $s, t \in (0, \tau_K]$, and $\Sigma(s, t)=E\{\eta(s)\eta(t)^{\rm T}\}$ with $\eta(t)=\bm Z\big\{N(\exp\{\bm Z^{\rm T}\bbeta_0(t)\})-\int_{0}^{t}I[X_{i}\ge \exp\{\bm Z^{\rm T}\bbeta_{0}(u)\}] {\rm d}H(u)\big\}$. \\

The proofs of Theorems 1 and 2 are sketched in Section B of the Supplementary Material.
Specifically, in contrast with    the proof of Theorem 2 in \cite{Peng2008},   the above asymptotic result implies that the proposed online   estimator $\widehat{\bbeta}_{N_B}(\tau)$ based
only on current data batch $D_B$ and some summary statistics
from historical data can achieve the same statistical efficiency as the oracle
estimator $\widetilde{\bbeta}_{N_B}(\tau)$, which is  calculated with all the  raw  data   up to the accumulation
point $B$.

\subsection{Variance estimation}
\label{VarianceEstimation}
To make inferences about $\bm \beta_0(\tau)$ in the  proposed online method, one often needs to estimate the covariance matrix of   $\widehat{\bbeta}_{N_B}(\tau)$.
However,    the asymptotic  covariance  matrix derived in Theorem 2 involves unknown density functions $f(\cdot \mid \bm Z)$ and
$\widetilde{f}(\cdot \mid \bm Z)$, making the traditional plug-in variance estimation procedure difficult to use in finite samples.
As an alternative, we herein  extend  the resampling approach used in  \cite{Peng2008} to the online setting.
In particular, define $$\widetilde{S}_{n}(\bbeta(\tau))=\frac{1}{n}\sum_{i=1}^{n}\bm Z_{i}\left(N_{i}(\exp\{\bm Z_{i}^{\rm T}\bbeta(\tau)\})-\int_{0}^{\tau}I\left[X_{i}\ge \exp\{\bm Z_{i}^{\rm T}\bbeta(u)\}\right]{\rm d}H(u)\right).$$
According to Theorem 2 and its proof given in Section B of the Supplementary Material, we have
$$
\sqrt{N_{B}}\{\widehat{\bbeta}_{N_B}(\tau)-\bbeta_0(\tau)\}=\Gamma_{N_B}^{-1} (\tau)\Phi \{-\sqrt{N_B} \widetilde{S}_{N_B}(\bbeta_0(\tau))\}+o_p(1).
$$
Because $\Phi$ is a linear operator, we have
$$
\widehat{\bbeta}_{N_B}(\tau)=\frac{\bbeta_0(\tau)}{\sqrt{N_{B}}}+\Gamma_{N_B}^{-1} (\tau)\frac{1}{N_B}\sum_{b=1}^{B}\Phi \{-n_b\widetilde{S}_{n_b}(\bbeta_0(\tau))\}+o_{p}\Big(\frac{1}{\sqrt{N_{B}}}\Big).
$$
Since $\sqrt{n_b}\widetilde{S}_{n_b}(\widetilde{\bbeta}_{n_b}(\tau))=o_p(1)$, we have
$$
\begin{aligned}
\widehat{\bbeta}_{N_B}(\tau)=&\frac{\bbeta_0(\tau)}{\sqrt{N_{B}}}+\Gamma_{N_B}^{-1} (\tau)\frac{1}{N_B}\sum_{b=1}^{B}\Phi \{n_b\widetilde{S}_{n_b}(\widetilde{\bbeta}_{n_b}(\tau))-n_b\widetilde{S}_{n_b}(\bbeta_0(\tau))\}\\
&+o_{p}\Big(\frac{1}{\sqrt{N_{B}}}\Big)+o_p\Big(\frac{1}{N_B}\sum_{b=1}^{B}\sqrt{n_b}\Big)\\
=&\frac{\bbeta_0(\tau)}{\sqrt{N_{B}}}+\Gamma_{N_B}^{-1} (\tau)\frac{1}{N_B}\sum_{b=1}^{B}n_b\Gamma_{n_b}(\tau)\{\widetilde{\bbeta}_{n_b}(\tau)-\bbeta_0(\tau)\}\\
&+o_{p}\Big(\frac{1}{\sqrt{N_{B}}}\Big)+o_p\Big(\frac{1}{N_B}\sum_{b=1}^{B}\sqrt{n_b}\Big).
\end{aligned}
$$
Let $V_{n_b}(\tau)$ be the covariance matrix of $\widetilde{\bbeta}_{n_b}(\tau)$  based on the $b$th data block,  the covariance matrix of $\widehat{\bbeta}_{N_B}(\tau)$ is given by
$$
\begin{aligned}
Var(\widehat{\bbeta}_{N_B}(\tau)) =&\frac{1}{N_B^2}\Gamma_{N_B}(\tau)^{-1} \Bigg[\sum_{b=1}^{B}n_b^2 \Gamma_{n_b}(\tau)V_{n_b}(\tau)\Gamma_{n_b}(\tau)\Bigg] \Gamma_{N_B}(\tau)^{-1}\\
  =& \widetilde{\Gamma}_{N_B}(\tau)^{-1} \Bigg[\sum_{b=1}^{B}\widetilde{\Gamma}_{n_b}(\tau)V_{n_b}(\tau)\widetilde{\Gamma}_{n_b}(\tau)\Bigg] \widetilde{\Gamma}_{N_B}(\tau)^{-1}\\
 =&\Big\{\widetilde{\Gamma}_{N_{B-1}}(\tau)+\widetilde{\Gamma}_{n_{B}}(\tau)\Big\}^{-1} \Bigg[\sum_{b=1}^{B-1}\widetilde{\Gamma}_{n_b}(\tau)V_{n_b}(\tau)\widetilde{\Gamma}_{n_b}(\tau)+\widetilde{\Gamma}_{n_B}(\tau)V_{n_B}(\tau)\widetilde{\Gamma}_{n_B}(\tau)\Bigg] \\
 & \times \Big\{\widetilde{\Gamma}_{N_{B-1}}(\tau)+\widetilde{\Gamma}_{n_{B}}(\tau)\Big\}^{-1},
\end{aligned}
$$
where $\widetilde{\Gamma}_{N_B}(\tau)=N_B \Gamma_{N_B}(\tau)$ and $\widetilde{\Gamma}_{n_b}(\tau)=n_b \Gamma_{n_b}(\tau)$.

Thus, to calculate $Var(\widehat{\bbeta}_{N_B}(\tau))$
at the accumulation point $B$, we need to store $\widetilde{\Gamma}_{N_{B-1}}(\tau)$ and $\sum_{b=1}^{B-1}\widetilde{\Gamma}_{n_b}(\tau)V_{n_b}(\tau)\widetilde{\Gamma}_{n_b}(\tau)$ from previous calculations, and compute $\widetilde{\Gamma}_{n_B}(\tau)$ and $\widetilde{\Gamma}_{n_B}(\tau) $ $\times V_{n_B}(\tau) \widetilde{\Gamma}_{n_B}(\tau)$ with the current data block, where the estimates of $\Gamma_{N_B}(\tau)$  and $\Gamma_{n_b}(\tau)$  are given in Section \ref{ReLS}.
For each $b= 1, \ldots, B$,
since $V_{n_b}(\tau)$ consists of the unknown density functions $f( \cdot \mid \bm Z)$ and $\widetilde{f}(\cdot \mid \bm Z)$, we propose to estimate it by a resampling approach based on the $b$th  data batch.
More specifically,
let $\widetilde{w}_1, \ldots, \widetilde{w}_{n_b}$ be a set of independent nonnegative weights generated from a   parametric distribution with mean 1 and variance 1. For the $b$th   data batch with $b= 1, \ldots, B$, we define
$$
\begin{aligned}
\bar{l}_{kn_b}(\bbeta(\tau_k))=& \frac{1}{n_b}\sum_{i=1}^{n_b}\Bigg[\Delta_{i}|\widetilde{w}_i \log X_{i}-\widetilde{w}_i \bbeta(\tau_k)^{\rm T}\bm Z_{i}|+\bbeta(\tau_k)^{\rm T}\widetilde{w}_i\Delta_{i} \bm  Z_{i}\Bigg.   \\ &
\Bigg.-2\bbeta(\tau_k)^{\rm T}\widetilde{w}_i\bm  Z_{i}\Bigg(\sum_{r=0}^{k-1}I\Big[ X_{i}\ge \exp\{\bm  Z_{i}^{\rm T}\widetilde{\bbeta}(\tau_r)\}\Big]\{H(\tau_{r+1})-H(\tau_{r})\}\Bigg)\Bigg].
\end{aligned}
$$
By minimizing $\overline{l}_{kn_b}(\bbeta(\tau_k))$, one can obtain an estimator of $\bbeta(\tau_k)$ defined by $\bbeta_{n_b}^*(\tau_k)$ for each $k=1, \ldots, K$.
As in Section \ref{Review},
 $\bbeta_{n_b}^*(\tau)$ with $\tau \in (0, \tau_K]$ can be regarded  as a right-continuous step function that jumps only at $\tau_k$ with $k=1, \ldots, K$.
 For each $b= 1, \ldots, B$,
through repeatedly generating the weights $\{\widetilde{w}_1, \ldots, \widetilde{w}_{n_b}\}$ for $S$ times, we can obtain $S$ realizations of  $\bbeta_{n_b}^{*} (\tau)$, denoted by $\{ \bbeta_{n_b}^{*s}(\tau)\}_{s=1}^{S}$.
Thus, for each $b = 1, \ldots, B$,
the covariance matrix of $\widetilde{\bbeta}_{n_b}(\tau)$ denoted by $Var(\widehat{\bbeta}_{n_b}(\tau))$  for a specific     $\tau \in (0, \tau_K]$ can thus be estimated by the sample covariance matrix of $\{\bbeta_{n_b}^{*s}(\tau)\}_{s=1}^{S}$.

\section{Simulation studies}
We conducted extensive simulations to evaluate the empirical
performance of the proposed online method.
We also  processed the entire simulated raw data once with
the oracle  method of  \cite{Peng2008}, which can essentially serve as a benchmark for comparison regarding the estimation accuracy and efficiency.  Let the batch size $B=5, 20$ or $40$, and the sample size of  the $b$th data batch be $n_b=1000$   for $b=1, \ldots, B$.
In the first simulation,  for each batch  $b$, we generated the event times
from an accelerated failure time (AFT) model with i.i.d. errors,
$$\log T = \gamma_1 \widetilde{Z}_1  +  \gamma_2  \widetilde{Z}_2   + \epsilon,$$
where $\widetilde{Z}_1 \sim Uniform(0,1)$,  $\widetilde{Z}_2 \sim Bernoulli(0.5)$,
$(\gamma_1,  \gamma_2) = (0.5,-0.5)$
and  the error  $\epsilon$ follows either the  standard normal distribution $N(0,1)$
or the skewed  extreme value distribution with scale  1, denoted by $Extreme (1)$.
In this setting with $\bZ = (1,  \widetilde{Z}_1,   \widetilde{Z}_2)^{\rm T}$,   the underlying regression coefficient
is  $\bbeta (\tau) = (\beta^{(1)} (\tau),
\beta^{(2)} (\tau), \beta^{(3)} (\tau))^{\rm T} = (Q_{\epsilon}(\tau),  \gamma_1,  \gamma_2)^{\rm T}$, where $Q_{\epsilon}(\tau)$ denotes
the $\tau$th quantile of $\epsilon$.
Notably, in this AFT model, $\bbeta^{(1)} (\tau)$ varies with $\tau$  whereas $\bbeta^{(2)} (\tau)$ and $\bbeta^{(3)}(\tau)$ are constants.
The  censoring time $C$ follows   $Uniform(0, U_c)$, where $U_c$  was selected to produce  about $50\%$ right censoring rate.   The simulation results presented below are based on  $500$  replications.

In the proposed online method,
we estimated $\Gamma_{n_b}(\tau)$ and $V_{n_b}(\tau)$   given in Sections  \ref{ReLS}  and  \ref{VarianceEstimation}  by the ReLS and  resampling methods with 250 replications, respectively.
In particular, we independently generated $\bm \xi_k$  from $N(0,1)$ and the resampling weights   from   $Exponential(1)$  to estimate   ${\Gamma}_{n_b}(\tau)$ and $V_{n_b}(\tau)$, respectively.
For $b=2, \ldots, B$, we used the estimator  obtained from the online estimation at  the  accumulation point $(b-1)$   as the initial value  to implement the MM algorithm for each quantile level in the  online estimation at  the  accumulation point $b$.
The quantile grid size was equally set to be $0.01$.
Tables 1 and 2 present the simulation results for $\tau=$ 0.1, 0.3 or 0.5 when the error term $\epsilon$ follows $N(0,1)$ and $Extreme (1)$, respectively.
The tables include the  empirical bias  (Bias) calculated by the average of 500 regression parameter estimates minus the true value, the standard deviation of the 500 regression parameter  estimates (SD), the average of 500 standard error  estimates (ASE), the $95\%$ empirical coverage probability based on the normal approximation (CP), and the computing time given in minutes for a replication (Time).


Tables 1 and 2 indicate  that the proposed online estimators are
virtually unbiased across all quantiles and batch sizes considered.
The ratio of SD and ASE is consistently close to 1, manifesting a reliable
performance of the proposed variance estimation procedure.
The coverage probabilities are all around the nominal level $95\%$.
As anticipated,  the proposed online estimators become more efficient as the batch size $B$ increases.
Moreover, it is worth noting that the proposed online method
achieves similar estimation accuracy as the oracle method, which is based on the entire raw data, but at a
slight expense of estimation efficiency.
Regarding the computation time,  compared to the  oracle method, the proposed online method
took  a longer time to implement  due to the use of two resampling-based algorithms when the batch size $B=5$ and 20,
whereas it exhibited faster computations when $B$ was increased to   40.
These comparative results manifest the advantages of the proposed online method (e.g.,  computation time and   required computer memory)    in handling  large-scale streaming data sets.

\begin{table}[htbp]
    \centering
    {\small 
    {Table 1: Simulation results with $\epsilon \sim N(0,1)$.}}\\
        \vspace{0.2cm}
\scalebox{0.75}{
    \begin{tabular}{lcccccccccccccccccccc}
    \hline
         &           &       &      &\multicolumn{4}{c}{$\beta^{(1)}$}       & & \multicolumn{4}{c}{$\beta^{(2)}$} & & \multicolumn{4}{c}{$\beta^{(3)}$} \\
        \cline{5-8} \cline{10-13} \cline{15-18}
$B$ &Method&Time& $\tau$& Bias & SD & ASE & CP & & Bias & SD &  ASE & CP & & Bias & SD &  ASE & CP\\ \hline
    5	&	Online	&	3.9 	&	0.1	&	0.008 	&	0.070 	&	0.067 	&	0.930 	&	&	0.005 	&	0.110 	&	0.105 	&	0.923 	&	&	-0.007 	&	0.062 	&	0.054 	&	0.930 	\\
	&		&		&	0.3	&	0.008 	&	0.051 	&	0.051 	&	0.933 	&	&	0.002 	&	0.086 	&	0.084 	&	0.923 	&	&	-0.002 	&	0.044 	&	0.043 	&	0.946 	\\
	&		&		&	0.5	&	0.015 	&	0.067 	&	0.063 	&	0.930 	&	&	-0.002 	&	0.098 	&	0.094 	&	0.927 	&	&	-0.006 	&	0.049 	&	0.048 	&	0.923 	\\
	&	Oracle	&	1.1 	&	0.1	&	0.008 	&	0.055 	&	0.057 	&	0.948 	&	&	0.004 	&	0.088 	&	0.089 	&	0.946 	&	&	-0.004 	&	0.052 	&	0.051 	&	0.930 	\\
	&		&		&	0.3	&	0.011 	&	0.045 	&	0.046 	&	0.950 	&	&	-0.001 	&	0.071 	&	0.072 	&	0.944 	&	&	-0.003 	&	0.043 	&	0.042 	&	0.952 	\\
	&		&		&	0.5	&	0.015 	&	0.047 	&	0.047 	&	0.932 	&	&	0.002 	&	0.076 	&	0.074 	&	0.928 	&	&	-0.004 	&	0.044 	&	0.043 	&	0.934 	\\
20	&	Online	&	10.9 	&	0.1	&	0.007 	&	0.033 	&	0.037 	&	0.960 	&	&	0.004 	&	0.054 	&	0.058 	&	0.960 	&	&	-0.002 	&	0.029 	&	0.029 	&	0.957 	\\
	&		&		&	0.3	&	0.006 	&	0.027 	&	0.026 	&	0.933 	&	&	0.002 	&	0.040 	&	0.042 	&	0.963 	&	&	-0.001 	&	0.021 	&	0.022 	&	0.953 	\\
	&		&		&	0.5	&	0.013 	&	0.031 	&	0.030 	&	0.931 	&	&	0.002 	&	0.044 	&	0.044 	&	0.947 	&	&	-0.004 	&	0.026 	&	0.024 	&	0.925 	\\
	&	Oracle	&	8.9 	&	0.1	&	0.007 	&	0.028 	&	0.028 	&	0.942 	&	&	0.003 	&	0.045 	&	0.044 	&	0.942 	&	&	-0.002 	&	0.027 	&	0.026 	&	0.938 	\\
	&		&		&	0.3	&	0.007 	&	0.021 	&	0.023 	&	0.952 	&	&	0.005 	&	0.035 	&	0.036 	&	0.950 	&	&	-0.002 	&	0.020 	&	0.021 	&	0.954 	\\
	&		&		&	0.5	&	0.011 	&	0.022 	&	0.024 	&	0.932 	&	&	0.002 	&	0.036 	&	0.037 	&	0.950 	&	&	-0.004 	&	0.020 	&	0.021 	&	0.954 	\\
40	&	Online	&	16.6 	&	0.1	&	0.005 	&	0.028 	&	0.029 	&	0.967 	&	&	0.003 	&	0.045 	&	0.045 	&	0.953 	&	&	-0.001 	&	0.019 	&	0.022 	&	0.970 	\\
	&		&		&	0.3	&	0.006 	&	0.019 	&	0.021 	&	0.956 	&	&	0.001 	&	0.031 	&	0.034 	&	0.967 	&	&	-0.001 	&	0.015 	&	0.017 	&	0.960 	\\
	&		&		&	0.5	&	0.012 	&	0.025 	&	0.024 	&	0.922 	&	&	0.002 	&	0.041 	&	0.038 	&	0.937 	&	&	-0.001 	&	0.019 	&	0.018 	&	0.950 	\\
	&	Oracle	&	24.1	&	0.1	&	0.006 	&	0.019 	&	0.019 	&	0.923 	&	&	0.003 	&	0.030 	&	0.031 	&	0.950 	&	&	-0.002 	&	0.018 	&	0.018 	&	0.933 	\\
	&		&		&	0.3	&	0.006 	&	0.016 	&	0.016 	&	0.916 	&	&	0.002 	&	0.025 	&	0.025 	&	0.943 	&	&	-0.001 	&	0.015 	&	0.015 	&	0.940 	\\
	&		&		&	0.5	&	0.012 	&	0.017 	&	0.017 	&	0.887 	&	&	0.002 	&	0.027 	&	0.026 	&	0.933 	&	&	-0.004 	&	0.015 	&	0.015 	&	0.936 	\\

         \hline
    \end{tabular}
    }
             \begin{center}
{\small Note:   Bias,     empirical bias;  SD,  standard deviation of the 500 regression parameter  estimates; ASE,   average of 500 standard error  estimates;  CP,  $95\%$ coverage probability; Time, computation time  for a replication given in minutes.
}
\end{center}
\end{table}

\begin{table}[htbp]
    \centering
    {\small 
    {Table 2: Simulation results with $\epsilon$ following the  Extreme value distribution.}}\\
        \vspace{0.2cm}
\scalebox{0.75}{
    \begin{tabular}{lcccccccccccccccccccc}
    \hline
         &           &       &      &\multicolumn{4}{c}{$\beta^{(1)}$}       & & \multicolumn{4}{c}{$\beta^{(2)}$} & & \multicolumn{4}{c}{$\beta^{(3)}$} \\
        \cline{5-8} \cline{10-13} \cline{15-18}
$B$ &Method&Time& $\tau$& Bias & SD & ASE & CP & & Bias & SD &  ASE & CP & & Bias & SD &  ASE & CP\\ \hline
    5	&	Online	&	3.6 	&	0.1	&	0.005 	&	0.050 	&	0.049 	&	0.947 	&	&	0.007 	&	0.080 	&	0.079 	&	0.930 	&	&	-0.003 	&	0.042 	&	0.041 	&	0.930 	\\
	&		&		&	0.3	&	0.010 	&	0.051 	&	0.049 	&	0.943 	&	&	0.008 	&	0.077 	&	0.079 	&	0.940 	&	&	-0.002 	&	0.044 	&	0.041 	&	0.937 	\\
	&		&		&	0.5	&	0.015 	&	0.067 	&	0.063 	&	0.943 	&	&	0.007 	&	0.103 	&	0.108 	&	0.943 	&	&	-0.008 	&	0.057 	&	0.052 	&	0.921 	\\
	&	Oracle	&	1.0 	&	0.1	&	0.007 	&	0.045 	&	0.044 	&	0.927 	&	&	0.003 	&	0.067 	&	0.068 	&	0.937 	&	&	-0.004 	&	0.037 	&	0.039 	&	0.953 	\\
	&		&		&	0.3	&	0.009 	&	0.044 	&	0.044 	&	0.950 	&	&	0.006 	&	0.068 	&	0.069 	&	0.943 	&	&	-0.006 	&	0.040 	&	0.040 	&	0.947 	\\
	&		&		&	0.5	&	0.015 	&	0.055 	&	0.054 	&	0.923 	&	&	0.011 	&	0.082 	&	0.085 	&	0.943 	&	&	-0.007 	&	0.049 	&	0.048 	&	0.937 	\\
20	&	Online	&	11.9 	&	0.1	&	0.004 	&	0.026 	&	0.026 	&	0.930 	&	&	0.005 	&	0.040 	&	0.040 	&	0.940 	&	&	-0.002 	&	0.020 	&	0.021 	&	0.960 	\\
	&		&		&	0.3	&	0.008 	&	0.025 	&	0.025 	&	0.927 	&	&	0.001 	&	0.040 	&	0.039 	&	0.947 	&	&	-0.003 	&	0.021 	&	0.021 	&	0.953 	\\
	&		&		&	0.5	&	0.013 	&	0.039 	&	0.035 	&	0.927 	&	&	0.004 	&	0.062 	&	0.055 	&	0.940 	&	&	-0.005 	&	0.031 	&	0.028 	&	0.927 	\\
	&	Oracle	&	10.0 	&	0.1	&	0.005 	&	0.021 	&	0.022 	&	0.950 	&	&	0.002 	&	0.032 	&	0.034 	&	0.967 	&	&	-0.003 	&	0.019 	&	0.020 	&	0.947 	\\
	&		&		&	0.3	&	0.008 	&	0.022 	&	0.022 	&	0.943 	&	&	0.003 	&	0.035 	&	0.034 	&	0.927 	&	&	-0.003 	&	0.019 	&	0.020 	&	0.940 	\\
	&		&		&	0.5	&	0.014 	&	0.028 	&	0.027 	&	0.923 	&	&	0.003 	&	0.043 	&	0.042 	&	0.943 	&	&	-0.003 	&	0.025 	&	0.024 	&	0.927 	\\
40	&	Online	&	16.7 	&	0.1	&	0.004 	&	0.021 	&	0.021 	&	0.953 	&	&	0.002 	&	0.031 	&	0.032 	&	0.963 	&	&	-0.002 	&	0.015 	&	0.016 	&	0.967 	\\
	&		&		&	0.3	&	0.007 	&	0.018 	&	0.020 	&	0.940 	&	&	0.002 	&	0.027 	&	0.031 	&	0.972 	&	&	-0.003 	&	0.014 	&	0.014 	&	0.950 	\\
	&		&		&	0.5	&	0.012 	&	0.027 	&	0.027 	&	0.946 	&	&	0.001 	&	0.044 	&	0.042 	&	0.940 	&	&	-0.004 	&	0.022 	&	0.021 	&	0.950 	\\
	&	Oracle	&	21.1	&	0.1	&	0.004 	&	0.016 	&	0.015 	&	0.920 	&	&	0.002 	&	0.025 	&	0.024  	&	0.930 	&	&	-0.002 	&	0.013 	&	0.014 	&	0.956 	\\
	&		&		&	0.3	&	0.007 	&	0.017 	&	0.016 	&	0.890 	&	&	0.003 	&	0.025 	&	0.024 	&	0.936 	&	&	-0.002 	&	0.013 	&	0.014 	&	0.950 	\\
	&		&		&	0.5	&	0.010 	&	0.019 	&	0.019 	&	0.900 	&	&	0.003 	&	0.030 	&	0.030 	&	0.940 	&	&	-0.003 	&	0.016 	&	0.017 	&	0.967 	\\

         \hline
    \end{tabular}
    }
             \begin{center}
{\small Note:   Bias,     empirical bias;  SD,  standard deviation of the 500 regression parameter  estimates; ASE,   average of 500 standard error  estimates;  CP,  $95\%$ coverage probability; Time, computation time for a replication given in minutes.
}
\end{center}
\end{table}

In the second   simulation,
we investigated a setting where the event times were generated from a regression model
with i.i.d. heteroscedastic errors.
Specifically, we considered the model
$\log T =\gamma_1  \widetilde{Z}_1 + \gamma_2 \widetilde{Z}_2 \zeta+\epsilon$,
where    $\widetilde{Z}_1 \sim Uniform(0,1)$,  $\widetilde{Z}_2 \sim Bernoulli(0.5)$,
$(\gamma_1,  \gamma_2) = (0.5,-0.5)$,
 $\epsilon$ follows  the  standard normal distribution $N(0, 1)$  and
 $\zeta$ follows the exponential  distribution with mean 1 denoted by $Exponential (1)$.
In the setting, $\bbeta^{(1)}(\tau)$ and $\bbeta^{(3)}(\tau)$ in the quantile regression model \eqref{m1} depend on $\tau$ and  $\bbeta^{(2)}(\tau)$ equals the constant $0.5$.
We  generated the  censoring time $C$ from  $Uniform(0, 2.5)$,   leading to a censoring rate of about $50\%$.
The results summarized in Table 3 indicate that the  proposed online method in this heteroscedastic error case
still performs reasonably and  is comparable to the oracle method in terms of the estimation accuracy and efficiency.
Like the first simulation,
the proposed online method exhibited an advantage  regarding the computation speed when the number of data batches $B$ is relatively large (e.g., 40).

\begin{table}[htbp]
    \centering
    {\small 
    {Table 3: Simulation results with the heteroscedastic $\epsilon$.}}\\
        \vspace{0.2cm}
\scalebox{0.75}{
    \begin{tabular}{lcccccccccccccccccccc}
    \hline
         &           &       &      &\multicolumn{4}{c}{$\beta^{(1)}$}       & & \multicolumn{4}{c}{$\beta^{(2)}$} & & \multicolumn{4}{c}{$\beta^{(3)}$} \\
        \cline{5-8} \cline{10-13} \cline{15-18}
$B$ &Method&Time& $\tau$& Bias & SD & ASE & CP & & Bias & SD &  ASE & CP & & Bias & SD &  ASE & CP\\ \hline
   5	&	Online	&	3.6	&	0.1	&	0.008 	&	0.073 	&	0.072 	&	0.931 	&	&	0.004 	&	0.124 	&	0.116 	&	0.924 	&	&	-0.002 	&	0.061 	&	0.060 	&	0.927 	\\
	&		&		&	0.3	&	0.008 	&	0.051 	&	0.053 	&	0.962 	&	&	0.007 	&	0.085 	&	0.086 	&	0.938 	&	&	-0.003 	&	0.045 	&	0.045 	&	0.948 	\\
	&		&		&	0.5	&	0.013 	&	0.061 	&	0.058 	&	0.924 	&	&	0.013 	&	0.099 	&	0.094 	&	0.927 	&	&	-0.003 	&	0.049 	&	0.048 	&	0.948 	\\
	&	Oracle	&	1.4	&	0.1	&	0.011 	&	0.061 	&	0.059 	&	0.926 	&	&	0.002 	&	0.095 	&	0.097 	&	0.950 	&	&	-0.010 	&	0.057 	&	0.056 	&	0.940 	\\
	&		&		&	0.3	&	0.007 	&	0.048 	&	0.048 	&	0.954 	&	&	0.006 	&	0.077 	&	0.076 	&	0.940 	&	&	-0.005 	&	0.042 	&	0.044 	&	0.940 	\\
	&		&		&	0.5	&	0.010 	&	0.046 	&	0.049 	&	0.950 	&	&	0.007 	&	0.072 	&	0.078 	&	0.960 	&	&	-0.005 	&	0.046 	&	0.045 	&	0.926 	\\
20	&	Online	&	9.3	&	0.1	&	0.007 	&	0.037 	&	0.038 	&	0.941 	&	&	0.001 	&	0.061 	&	0.062 	&	0.952 	&	&	-0.006 	&	0.028 	&	0.032 	&	0.960 	\\
	&		&		&	0.3	&	0.008 	&	0.026 	&	0.028 	&	0.958 	&	&	0.002 	&	0.043 	&	0.045 	&	0.952 	&	&	-0.003 	&	0.024 	&	0.024 	&	0.939 	\\
	&		&		&	0.5	&	0.011 	&	0.030 	&	0.030 	&	0.923 	&	&	0.002 	&	0.051 	&	0.048 	&	0.935 	&	&	-0.003 	&	0.027 	&	0.025 	&	0.924 	\\
	&	Oracle	&	7.6	&	0.1	&	0.007 	&	0.029 	&	0.030 	&	0.960 	&	&	0.002 	&	0.045 	&	0.048 	&	0.954 	&	&	-0.006 	&	0.026 	&	0.028 	&	0.952 	\\
	&		&		&	0.3	&	0.008 	&	0.024 	&	0.023 	&	0.940 	&	&	0.003 	&	0.036 	&	0.037 	&	0.954 	&	&	-0.003 	&	0.022 	&	0.022 	&	0.942 	\\
	&		&		&	0.5	&	0.009 	&	0.024 	&	0.024 	&	0.923 	&	&	0.002 	&	0.041 	&	0.038 	&	0.920 	&	&	-0.004 	&	0.023 	&	0.022 	&	0.944 	\\
40	&	Online	&	13.9	&	0.1	&	0.005 	&	0.029 	&	0.031 	&	0.966 	&	&	0.002 	&	0.048 	&	0.050 	&	0.970 	&	&	-0.005 	&	0.022 	&	0.022 	&	0.962 	\\
	&		&		&	0.3	&	0.007 	&	0.019 	&	0.022 	&	0.966 	&	&	0.002 	&	0.032 	&	0.036 	&	0.968 	&	&	-0.002 	&	0.016 	&	0.017 	&	0.956 	\\
	&		&		&	0.5	&	0.008 	&	0.026 	&	0.024 	&	0.945 	&	&	-0.001 	&	0.038 	&	0.039 	&	0.951 	&	&	-0.002 	&	0.020 	&	0.019 	&	0.943 	\\
	&	Oracle	&	23.8	&	0.1	&	0.007 	&	0.020 	&	0.021 	&	0.936 	&	&	0.002 	&	0.032 	&	0.034 	&	0.960 	&	&	-0.007 	&	0.020 	&	0.020 	&	0.926 	\\
	&		&		&	0.3	&	0.008 	&	0.017 	&	0.017 	&	0.932 	&	&	0.002 	&	0.027 	&	0.027 	&	0.944 	&	&	-0.002 	&	0.015 	&	0.015 	&	0.940 	\\
	&		&		&	0.5	&	0.011 	&	0.017 	&	0.017 	&	0.910 	&	&	0.005 	&	0.028 	&	0.027 	&	0.952 	&	&	-0.002 	&	0.015 	&	0.016 	&	0.954 	\\

         \hline
    \end{tabular}
    }
             \begin{center}
{\small Note:   Bias,     empirical bias;  SD,  standard deviation of the 500 regression parameter  estimates; ASE,   average of 500 standard error  estimates;  CP,  $95\%$ coverage probability; Time, computation time  for a replication given in minutes.
}
\end{center}
\end{table}

\section{Application}

Launched in 1973,
the Surveillance, Epidemiology and End Results (SEER) Program of the National Cancer
Institute collected individual-level data on demographics, death and tumor-related information
of cancer cases in the states of Connecticut, Iowa, New Mexico, Utah and Hawaii, and the metropolitan areas of Detroit and San Francisco-Oakland.
For 50 years,  SEER  Program has continuously expanded to numerous additional  cancer registries and covered nearly 50\% of the U.S. population.
One can refer to \cite{FriedmanNegoita2024}
 for more details about SEER.
We applied our proposed online method to a colon cancer data set collected between 2006 and 2021 from  SEER program, which was available at
{\sl https:\scalebox{0.75}[0.8]{//}seer.cancer.gov\scalebox{0.75}[0.8]{/}data\scalebox{0.75}[0.8]{/}}.
The primary goal of our analysis was to investigate the potential   factors for the death risk of colon cancer patients.
Thus, the event time of interest $T$ was defined as the gap time from the  diagnosis of colon cancer to death.
If an individual was still alive at the last trial contact,  his/her death time was considered as right censoring.
This  colon cancer data set includes $450, 761$ observations  after deleting individuals with missing covariates, and the censoring rate is about $50.2\%$.

The covariates considered in our analysis include \texttt{Surgery}  (surgery treatment indicator, 1 for receiving a cancer-directed surgery by removing or
destroying the cancer tissue in the primary site  and 0 otherwise),
\texttt{Age} (age at the  diagnosis of colon cancer) and \texttt{Gender} (1 for male and 0 for female).
The continuous covariate   \texttt{Age} was  standardized to have mean 0 and variance 1.
Since the censoring rate is relatively high, we set $\tau_K = 0.5$ in view of the model identifiability.
The quantile grid size was equally set to be $0.01$.
In this application, streaming data were formed by yearly enrollment    from January 2006 to December 2021, with $B =16$ batches  and a total sample
size $N_B =450, 761$.
The sample size of each batch data varies from 24, 690 to 30, 287.
For a further comparison,  we also analyzed the  entire raw data with the oracle method  \citep{Peng2008}.
The analysis  results were showed in Figure 1, which displayed the estimated covariate effects with quantiles varying from $0.1$ to $ 0.5$  and  the 95\%
confidence bands.

\begin{figure}[h]
\centering
\includegraphics[height=5 in]{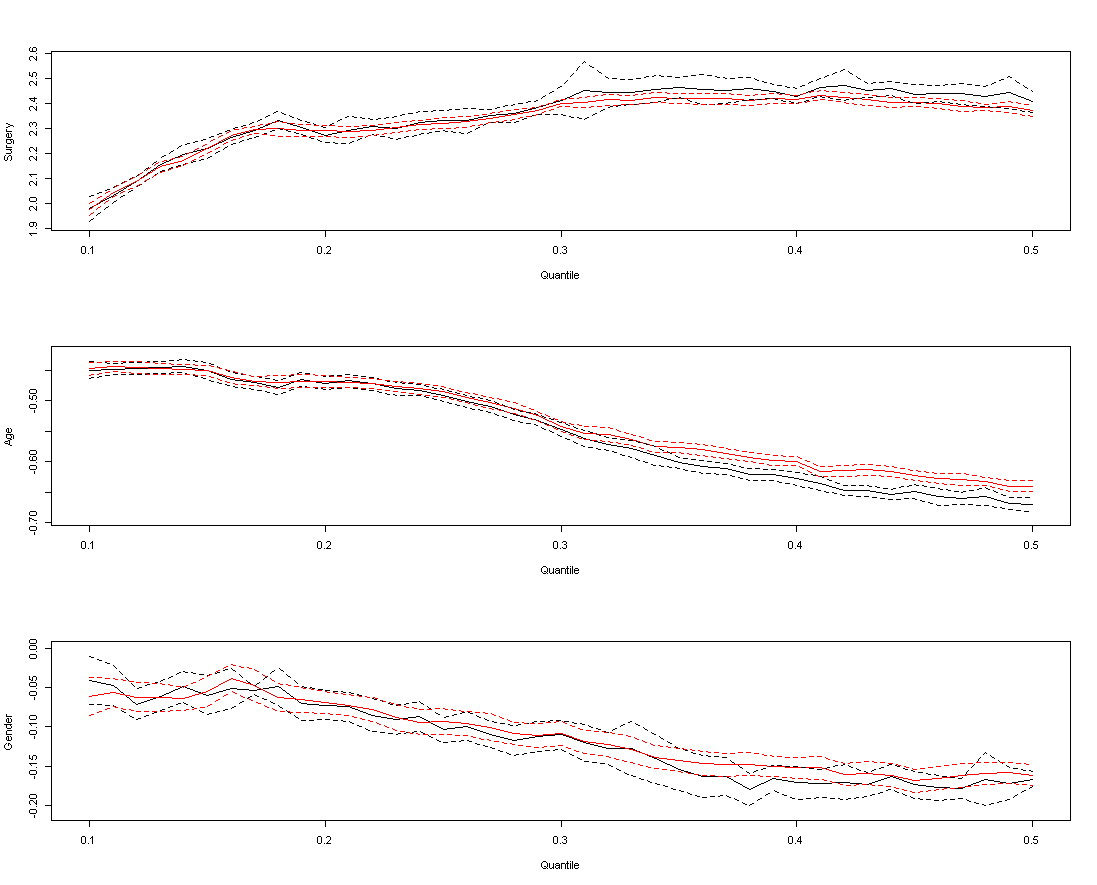}
\begin{center}
\caption{The analysis results of the colon cancer data, including the estimated covariate effects  for 0.1 to 0.5 quantiles (solid lines), and  95\%
confidence bands (dashed lines). The proposed   and oracle estimators are depicted by  black and red,
respectively.}\label{Realdata_optimal}
\end{center}
\end{figure}

%

Figure 1 shows that the covariate effect estimates and obtained inference conclusions  of our proposed online method are consistent  with those of  the oracle method, validating   our online method that only used the current data and
 summary statistics of historical data.
In particular,  \texttt{Surgery} exhibits    significant effects on the death risk at all quantile levels,  and  patients  who received surgery has   longer survival times or   lower death risk compared to those without surgery.
The effect trend of \texttt{Surgery} appears to be increasing as the quantile level increases, which  indicates that receiving surgery may  affect the survival experience of long term survivors more obviously but may have decreased impact on that of short term survivors.
 \texttt{Age} and  \texttt{Gender} are also  recognized as  significant,  which is  consistent with the conclusion given in \cite{2007Liang} and \cite{2013Purim}, and
 their
 effect trends are  downward as the quantile level increases.
 Such trends
 imply that  these two risk factors have more pronounced influences for  patients with shorter survival times.


In addition, it is worth noting that our proposed online  method takes 5.8 hours to analyze the colon cancer data, whereas   the oracle method requires 51.1 hours.
This attractive computational advantage arises because, in contrast with the oracle method, our proposed online method   only used the current data and
 summary statistics of historical data and thus release the computer memory.
Overall, through the above comparative analysis results,  one can clearly find that  the proposed online method remains statistically and computationally efficient when handing large-scale streaming data sets.

\section{Discussion}


In light of the appealing advantages of censored quantile regression,
we developed an online updating method   in the context of  streaming data sets.
The proposed method exhibited appealing features in terms of estimation accuracy and computational efficiency  due to the
low storage requirements.
Theoretical  developments and computational complex analysis also validated our findings in the numerical experiments.
Since the proposed online method is the first attempt to investigate  censored quantile regression with streaming
data sets,
one can easily envision  some potential future research directions.
First, note that the proposed renewable estimation method was developed  under the assumption that all data batches arise from a same regression model. However, this assumption may not hold in some practical applications, such as the USA airline study \citep{2024Chen}.
Thus, it may be worthwhile to   extend our current  work to the setting of heterogeneous streaming data sets.

Second,  note that the  proposed online method   relies on the conditionally independent censoring assumption.
Nevertheless, such an  assumption  may be  problematic in some practical applications   where the censoring  is driven by some outcome-related reasons \citep{LiPeng2015}, and
simply ignoring the dependent censoring in the analysis may lead to biased estimators.
Thus, another potential research direction is to develop an online method that can account for dependent censoring times.
Third, streaming data sets may also accompany by   high dimensional covariates, especially in   the genomic  and neuroimaging studies \citep{2017Wang, 2020Zhang}.
Hence,
developing an inference for high dimensional censored quantile regression with streaming survival data deserves further explorations.
Finally, developing a fast online model checking procedure for   censored quantile regression  warrants further investigations since traditional test  procedures are   inapplicable when survival data arrive sequentially in chunks.\\

\section*{Acknowledgements.}
Shuwei Li's research was partially supported by  the National Nature Science Foundation of China (Grant No. 12471251).
Liuquan Sun's research was partially supported by the National Natural Science Foundation of China (Grant Nos. 12171463 and 12426673).
Baoxue Zhang's research was partially supported by the National Natural Science Foundation of China (Grant No. 12271370).



\end{document}